\documentclass[12pt]{article}
\usepackage{amsmath}
\usepackage{amssymb}
\usepackage{amstext}
\usepackage{amscd}

\newtheorem{teo}{Theorem}[section]

\newtheorem{conj}[teo]{Conjecture}

\newtheorem{defini}[teo]{Definition}

\newcommand{\Hom}{\mbox{Hom}}

\newcommand{\End}{{\rm End}}
\newcommand{\GSp}{{\rm GSp}}
\newcommand{\GL}{{\rm GL}}
\newcommand{\SL}{{\rm SL}}
\newcommand{\Sh}{{\rm Sh}}

\newcommand{\Aut}{{\rm Aut}}
\newcommand{\Gal}{{\rm Gal}}

\newcommand{\MT}{{\rm MT}}

\newcommand{\ab}{{\rm ab}}

\newcommand{\ad}{{\rm ad}}

\newcommand{\CC}{{\mathbb C}}
\newcommand{\RR}{{\mathbb R}}
\newcommand{\ZZ}{{\mathbb Z}}
\newcommand{\QQ}{{\mathbb Q}}
\newcommand{\NN}{{\mathbb N}}

\newcommand{\HH}{{\mathbb H}}

\newcommand{\GG}{{\mathbb G}}
\newcommand{\SSS}{{\mathbb S}}
\newcommand{\AAA}{{\mathbb A}}
\newcommand{\UUU}{{\mathbb U}}

\newcommand{\Hol}{{\rm Hol}}

\newcommand{\cF}{{\cal F}}

\newcommand{\OOO}{{\cal O}}

\newcommand{\wt}{\widetilde}

\newenvironment{prf}[1]{\trivlist
\item[\hskip \labelsep{\it
#1.\hspace*{.3em}}]}{~\hspace{\fill}~$\square$\endtrivlist}

\title{The Andr\'e-Oort conjecture via o-minimality}

\author{Christopher Daw
\footnote{University
College London, Department of Mathematics, Gower Street, WC1E 6BT London, United Kingdom, 
e-mail: c.daw@ucl.ac.uk}}
\date{\today}

\begin{document}
\maketitle

\section{Introduction}

{\it Shimura varieties} are a distinguished class of algebraic varieties that parameterise important objects from linear algebra called {\it Hodge structures}. Often these Hodge structures correspond to families of so-called {\it Abelian varieties}. 

Additional structure on a Shimura variety $S$ arises through the existence of certain {\it algebraic correspondences} on $S$, i.e. subvarieties of $S\times S$, called {\it Hecke correspondences}. We can think of these as one-to-many maps 
\begin{align*}
T:S\rightarrow S.
\end{align*}
We endow $S$ with a set of so-called {\it special subvarieties}, defined as the set of all connected components of {\it Shimura subvarieties} and the irreducible components of their images under Hecke correspondences. This is analogous to the case of Abelian varieties (respectively {\it algebraic tori}), where special subvarieties are the translates of {\it Abelian subvarieties} (respectively {\it subtori}) by {\it torsion points}. A key property of special subvarieties is that connected components of their intersections are themselves special subvarieties. Thus, any subvariety $Y$ of $S$ is contained in a smallest special subvariety. If this happens to be a connected component of $S$ itself, then we say that $Y$ is {\it Hodge generic} in $S$.

We refer to the special subvarieties of dimension zero as {\it special points}. Special subvarieties contain a Zariski (in fact, analytically) dense set of special points. The {\it Andr\'e-Oort conjecture} predicts that this property characterises special subvarieties:

\begin{conj}{\bf(Andr\'e-Oort)}
Let $S$ be a Shimura variety and let $\Sigma$ be a set of special points contained in $S$. Every irreducible component of the Zariski closure of $\cup_{s\in\Sigma}s$ in $S$ is a special subvariety.
\end{conj}

A connected component of $S$ arises as a quotient $\Gamma\backslash D$, where $D$ is a certain type of complex manifold called a {\it Hermitian symmetric domain}, and $\Gamma$ is a certain type of discrete subgroup of $\Hol(D)^+$ called a {\it congruence subgroup}. From now on, we will use $S$ to denote this component. 

By \cite{KUY}, \S3, there exists a semi-algebraic {\it fundamental domain} $\cF\subset D$ for the action of $\Gamma$. By \cite{KUY}, Theorem 1.2, when the uniformisation map
\begin{align*}
\pi:D\rightarrow S
\end{align*} 
is restricted to $\cF$, one obtains a function definable in the o-minimal structure $\RR_{\rm an,exp}$. Through these observations, the Andr\'e-Oort conjecture becomes amenable to tools from o-minimality.

The purpose of this article is to explain the so-called {\it Pila-Zannier strategy} for proving the Andr\'e-Oort conjecture. This strategy first arose in a proof of the {\it Manin-Mumford conjecture} \cite{PZ} and was first adapted to Shimura varieties by Pila \cite{P}. We will follow the outline given by Ullmo \cite{U} for $\mathcal{A}^r_6$, where $\mathcal{A}_g$ is the moduli space for {\it principally polarised} Abelian varieties of dimension $g$. 

The first step is to show that, if $Y$ is an irreducible, Hodge generic subvariety of $S$, then the union of all positive-dimensional, special subvarieties contained in $Y$ is not Zariski dense in $Y$. The second step is to show that all but finitely many special points in $Y$ lie on a positive-dimensional, special subvariety contained in $Y$. 

Both steps require the {\it hyperbolic Ax-Lindemann-Weierstrass conjecture}, a geometric statement itself amenable to proof via o-minimality. Other articles in these proceedings will explain this conjecture in detail along with its analogue in the case of an Abelian variety. Let us just mention that the conjecture was first proven in the cocompact case by Ullmo and Yafaev \cite{UY}, then by Pila and Tsimerman for $\mathcal{A}_g$ \cite{PT}, and finally by Klingler, Ullmo and Yafaev in the general case \cite{KUY}. 

Ullmo demonstrates the first step in his article \cite{U}. Therefore, the focus of this article will be the second step. The strategy will be to compare lower bounds for the size of Galois orbits of special points with upper bounds for the heights of their pre-images in the fundamental domain. One concludes by applying the {\it Pila-Wilkie counting theorem} \cite{PW}, which states that the number of algebraic points of degree at most $k$ and height at most $T$, in the complement of all connected, positive-dimensional, semi-algebraic subsets of a set $X$, definable in an o-minimal structure, is $\ll_{\epsilon,k,X} T^{\epsilon}$. 

First, however, we will provide a brief introduction to the theory of Shimura varieties, as formulated by Deligne in his foundational articles \cite{D} and \cite{D2}. Our introduction is not by any means intended to be a full treatment of the topic but rather a preparatory guide for graduate students approaching it for the first time. We refer the reader to \cite{M} for a comprehensive account of Shimura varieties and for further details regarding the topics introduced here. 

\section{Hermitian symmetric domains}

We are primarily interested in the connected components of Shimura varieties. These initially arise as quotients $\Gamma\backslash D$, where $D$ is a certain type of complex manifold called a {\it Hermitian symmetric domain}, and $\Gamma$ is a {\it congruence subgroup}, acting via holomorphic automorphisms. The protypical example is the case of the upper half-plane 
\begin{align*}
D=\HH:=\{z\in\CC:\Im(z)>0\}
\end{align*} 
and $\Gamma=\SL_2(\ZZ)$, where any element of $\SL_2(\RR)$ acts on $\HH$ by
\begin{align*}
\left(\begin{array}{cc} a & b\\ c & d\end{array}\right)\cdot z=\frac{az+b}{cz+d}.
\end{align*}

We refer the reader to \cite{M}, \S1 for a detailed introduction to Hermitian symmetric domains. We merely summarise the key points. Unfortunately, the definition is not particularly enlightening:
\begin{defini}
A Hermitian symmetric domain is a connected complex manifold $D$ such that
\begin{itemize}

\item $D$ is equipped with a Hermitian metric. 

\item The group $\Aut(D)$ of holomorphic isometries acts transitively on $D$.

\item There exists a point $\tau\in D$ and an involution $\varphi\in\Aut(D)$ such that $\tau$ is an isolated fixed point of $\varphi$.

\item $D$ is of non-compact type.

\end{itemize}
\end{defini}
For any topological group $G$, we deonte its {\it neutral component} by $G^+$. By this we mean the connected component of $G$ containing the identity element ${\rm id}\in G$. By \cite{M}, Lemma 1.5, $\Aut(D)^+$ acts transitively on $D$ and, by \cite{M}, Proposition 1.6, it coincides with $\Hol(D)^+$, where $\Hol(D)$ denotes the larger group of all holomorphic automorphisms. Note that, given the transitivity of the $\Aut(D)$ action, the third condition is true for all points $\tau\in D$.

Returning to our earlier example, 
\begin{align*}
\Hol(\HH)=\SL_2(\RR)/\{\pm\rm id\}.
\end{align*} 
Since $\SL_2(\RR)$ is connected, so is $\Hol(\HH)$ and it therefore coincides with $\Aut(\HH)$. The element
\begin{align*}
\varphi:=\left(\begin{array}{cc} 0 & 1\\ -1 & 0\end{array}\right)\in\SL_2(\RR)
\end{align*}
fixes only $i\in\HH$, whereas $\varphi^2=-\rm id$. Hence, the image of $\varphi$ in $\Aut(\HH)$ is an involution of $\HH$ with an isolated fixed point.

However, from the definition follows a key property of Hermitian symmetric domains: by \cite{M}, Theorem 1.9, if we denote by $\UUU(\RR)$ the {\it circle group} $\{z\in\CC:|z|=1\}$, then for each point $\tau\in D$ there exists a unique homomorphism 
\begin{align*}
u_{\tau}:\UUU(\RR)\rightarrow\Hol(D)^+
\end{align*}
such that, for all $z\in\UUU(\RR)$,
\begin{itemize}
\item $u_{\tau}(z)(\tau)=\tau$. 
\item $u_{\tau}(z)$ acts as multiplication by $z$ on the tangent plane of $D$ at $\tau$.  
\end{itemize}
For example, consider the point $i\in\HH$ and let
\begin{align*}
h_i:\UUU(\RR)\rightarrow\SL_2(\RR):z=a+ib\mapsto\left(\begin{array}{cc} a & b\\ -b & a\end{array}\right).
\end{align*}
Then, for all $z\in\UUU(\RR)$, $h_i(z)$ fixes $i$ and
\begin{align*}
\left.\frac{d}{dz}\left(\frac{az+b}{-bz+a}\right)\right|_i=\frac{a^2+b^2}{(a-bi)^2}=\frac{z}{\bar{z}}.
\end{align*}
Therefore, if we define 
\begin{align*}
u_i:\UUU(\RR)\rightarrow\SL_2(\RR)/\{\pm\rm id\}:z\mapsto h_i(\sqrt{z})\text{ mod}\pm\rm id,
\end{align*}  
which is well-defined since $h_i(-1)=-\rm id$, then $u_i(z)$ acts on the tangent plane of $\HH$ at $i$ as multiplication by $z$.

Furthermore, note that, if $g\in\Hol(D)^+$ and $\tau\in D$, then the uniqueness of $u_{g\tau}$ implies that it must be the conjugate
\begin{align*}
gu_{\tau}g^{-1}:z\mapsto gu_{\tau}(z)g^{-1}.
\end{align*}
Therefore, since $\Hol(D)^+$ acts transitively on $D$, if we fix a point $\tau_0\in D$, we have a $\Hol(D)^+$-equivariant bijection between $D$ and the $\Hol(D)^+$-conjugacy class of $u_{\tau_0}$.

\section{Conjugacy classes}

By \cite{M}, Proposition 1.7, for any Hermitian symmetric domain $D$, there exists a unique, adjoint, semisimple algebraic group $G$ over $\RR$ such that \begin{align*}
G(\RR)^+=\Hol(D)^+.
\end{align*} 
By a {\it linear algebraic group} $G$ over $\RR$, we simply mean a group that can be defined as a subgroup of $\GL_n(\RR)$ by real polynomials in the matrix coefficients. For example, $\UUU(\RR)$ is a linear algebraic group over $\RR$ whose elements may be realised as those 
\begin{align*}
\left(\begin{array}{cc} a & b\\ c & d\end{array}\right)\in\GL_2(\RR)
\end{align*}
such that $a=d$, $b=-c$ and $a^2+b^2=1$ (in particular, $\UUU(\RR)$ is contained in $\SL_2(\RR)$). However, since $\UUU(\RR)$ is defined by polynomials, we can think of $\UUU(\RR)$ as the {\it real points} of what is usually considered the algebraic group, which we denote $\UUU$. Then, for any $\RR$-algebra $A$, $\UUU(A)$ is simply the group of solutions in $A$ to the above polynomials. 

By a {\it semisimple} algebraic group we mean a connected (for the Zariski topology), linear algebraic group that is isogenous to a product of almost-simple subgroups. By a {\it simple} algebraic group we mean a connected, linear algebraic group that is not commutative and has no proper, normal, algebraic subgroups other than the identity. By an {\it almost-simple} subgroup we mean a subgroup that is a simple algebraic group modulo a finite centre. An {\it isogeny} between semisimple algebraic groups is a surjective morphism with finite kernel. Two semisimple algebraic groups $H_1$ and $H_2$ are called {\it isogenous} if there exist isogenies
\begin{align*}
H_1\leftarrow G\rightarrow H_2,
\end{align*}
for some semisimple algebraic group $G$. This is an equivalence relation.
By {\it adjoint} we are referring to a group with trivial centre and, for a linear algebraic group $G$, we write $G^{\ad}$ for $G$ modulo its centre.

As shown in \cite{M}, \S1, every representation 
\begin{align*}
\UUU(\RR)\rightarrow\GL_n(\RR)
\end{align*} 
is algebraic i.e. the image is given by polynomials in the matrix entries and can be written $\UUU\rightarrow\GL_n$. In particular, for any $\tau\in D$, we may consider the homomorphism 
\begin{align*}
u_{\tau}:\UUU(\RR)\rightarrow G(\RR)^+
\end{align*} 
as an algebraic morphism $u_{\tau}:\UUU\rightarrow G$, yielding a morphism 
\begin{align*}
u_{\tau}:\UUU(A)\rightarrow G(A)
\end{align*}
for any $\RR$-algebra $A$. 

The group $\UUU$ is connected, commutative and consists entirely of {\it semisimple elements}. By the latter condition we mean that, for any representation 
\begin{align*}
\UUU\rightarrow\GL_n,
\end{align*}
any element in the image of $\UUU(\CC)$ can be diagonalised by an element of $\GL_n(\CC)$. The fact that $\UUU$ is also commutative implies that the elements in the image of $\UUU(\CC)$ can be simultaneously diagonalised by a single element of $\GL_n(\CC)$. We refer to a linear algebraic group of this sort as a {\it torus}.

For any representation of $\UUU$, the eigenvalues are given by homomorphisms $\UUU_{\CC}\rightarrow\GG_m$ called {\it characters}, where we write $\UUU_{\CC}$ for $\UUU$ considered as an algebraic group over $\CC$ and $\GG_m$ for the algebraic group such that, for any $\CC$-algebra $A$, 
\begin{align*}
\GG_m(A)=A^{\times}:=\{a\in A:a{\rm\ is\ invertible\ in\ }A\}.
\end{align*}
The characters are algebraic since, by definition, they are one-dimensional representations. In this case, each character is of the form $z\mapsto z^n$, where $n\in\ZZ$. 

By \cite{M}, Theorem 1.21, the homomorphism $u_{\tau}$ always satisfies the following three properties:
\begin{itemize}
\item Only the characters $z\mapsto 1$, $z\mapsto z$ and $z\mapsto z^{-1}$ occur in the representation of $\UUU(\RR)$ on the Lie algebra $\mathfrak{g}_{\CC}$ of $G_{\CC}$.
\item Conjugation by $u_{\tau}(-1)$ is a Cartan involution of $G$.
\item $u_{\tau}(-1)$ maps to a non-trivial element in every simple factor of $G$.
\end{itemize}
The {\it Lie algebra} of $G_{\CC}$ is the tangent plane of $G(\CC)$ at the identity. One definition is the kernel of the map
\begin{align*}
G(\CC[\epsilon])\rightarrow G(\CC)
\end{align*}
induced by $\epsilon\mapsto 0$, where $\epsilon^2=1$. Then $G(\CC)$ acts on $\mathfrak{g}_{\CC}$ by conjugation. For the definition of a {\it Cartan involution} see \cite{M}, \S1. 

On the other hand, if $G$ is any adjoint, semisimple algebraic group over $\RR$ and $u:\UUU\rightarrow G$ is a homomorphism satisfying the above three properties, then the $G(\RR)^+$-conjugacy class of $u$ naturally has the structure of a Hermitian symmetric domain $D$, for which 
\begin{align*}
G(\RR)^+=\Hol(D)^+
\end{align*} 
and $u(-1)$ is the involution associated to $u$ when regarded as a point of $D$.

\section{The Deligne torus}

Let $\SSS$ denote the linear algebraic group over $\RR$ such that $\SSS(\RR)=\CC^{\times}$. Similar to the case of $\UUU$ we may realise the elements of $\SSS(\RR)$ as those
\begin{align*}
\left(\begin{array}{cc} a & b\\ c & d\end{array}\right)\in\GL_2(\RR)
\end{align*}
such that $a=d$, $b=-c$. This is also a torus, usually referred to as the {\it Deligne torus}, and we have a short exact sequence
\begin{align*}
1\rightarrow\GG_m\xrightarrow{w}\SSS\rightarrow\UUU\rightarrow 1,
\end{align*} 
which on real points corresponds to
\begin{align*}
1\rightarrow\RR^{\times}\xrightarrow{r\mapsto r^{-1}}\CC^{\times}\xrightarrow{z\mapsto z/\bar{z}}\UUU(\RR)\rightarrow 1.
\end{align*}
Therefore, any homomorphism $u:\UUU\rightarrow G$ yields a homomorphism \begin{align*}
h:\SSS\rightarrow G,
\end{align*}
defined by $h(z)=u(z/\bar{z})$. Furthermore, $\UUU(\RR)$ will act on $\mathfrak{g}_{\CC}$ via the characters $z\mapsto 1$, $z\mapsto z$ and $z\mapsto z^{-1}$ if and only if $\SSS(\RR)$ acts on $\mathfrak{g}_{\CC}$ via the characters $z\mapsto 1$, $z\mapsto z/\bar{z}$ and $z\mapsto\bar{z}/z$. 

Conversely, let $h:\SSS\rightarrow G$ be a homomorphism such that $\SSS$ acts on $\mathfrak{g}_{\CC}$ via the characters $z\mapsto 1$, $z\mapsto z/\bar{z}$ and $z\mapsto\bar{z}/z$. Then $w(\GG_m(\RR))$ acts trivially on $\mathfrak{g}_{\CC}$, which implies that $h$ is trivial on $w(\GG_m(\RR))$, since the adjoint representation of $G$ on $\mathfrak{g}$ is faithful. Thus, $h$ arises from a homomorphism $u:\UUU\rightarrow G$.

Therefore, to give a $G(\RR)^+$-conjugacy class $D$ of homomorphisms $u:\UUU\rightarrow G$ satisfying the above three properties is the same as to give a $G(\RR)^+$-conjugacy class $X^+$ of homomorphisms $h:\SSS\rightarrow G$ satisfying the following:
\begin{itemize}
\item Only the characters $z\mapsto 1$, $z\mapsto z/\bar{z}$ and $z\mapsto\bar{z}/z$ occur in the representation of $\SSS(\RR)$ on $\mathfrak{g}_{\CC}$.
\item Conjugation by $h(i)$ constitutes a Cartan involution of $G$.
\item The element $h(i)$ maps to a non-trivial element in every simple factor of $G$.
\end{itemize}

\section{Hodge structures}

Therefore, the question should be {\it why are we interested in such conjugacy classes of morphisms $h:\SSS\rightarrow G$?} To understand this, we require the notion of a Hodge structure. Below is a brief summary of the relevant definitions. For a more comprehensive account, we refer the reader to \cite{M}, \S2.

For a real vector space $V$, we define complex conjugation on 
\begin{align*}
V(\CC):=V\otimes_{\RR}{\CC}
\end{align*} 
by $\overline{v\otimes z}:=v\otimes\overline{z}$. A {\it Hodge decomposition} of $V$ is a decomposition
\begin{align*}
V(\CC)=\bigoplus_{(p,q)\in\ZZ\times\ZZ}V^{p,q}
\end{align*}
such that $\overline{V^{p,q}}=V^{q,p}$. A {\it Hodge structure} is a real vector space $V$ with a Hodge decomposition. The set of pairs $(p,q)$ such that $V^{p,q}\neq 0$ is called the {\it type} of the Hodge structure and we refer to a Hodge structure of type $(-1,0),(0,-1)$ as a {\it complex structure}. 

For each $n\in\ZZ$, 
\begin{align*}
\bigoplus_{p+q=n}V^{p,q}
\end{align*} 
is stable under complex conjugation and equal to $V_n(\CC)$ for some real subspace $V_n$ of $V$. The decomposition $V=\oplus_nV_n$ is called the {\it weight decomposition} of $V$. If $V=V_n$, then $V$ is said to have {\it weight $n$}. The {\it Hodge filtration} associated with a Hodge structure $V$ of weight $n$ is
\begin{align*}
F:=\{\cdots\supset F^p\supset F^{p+1}\supset\cdots\},\ F^p:=\oplus_{r\geq p}V^{r,n-r}.
\end{align*}

A $\ZZ$-(respectively $\QQ$-){\it Hodge structure} is a free $\ZZ$-module (respectively $\QQ$-vector space) $V$ of finite rank (respectively dimension) equipped with a Hodge decomposition of 
\begin{align*}
V(\RR):=V\otimes\RR
\end{align*}
such that the weight decomposition is defined over $\QQ$.

Recall that we can identify $\SSS$ with a closed subgroup of $\GL_2$ as follows: for any $\RR$-algebra $A$, we realise $\SSS(A)$ as those matrices of the form
\begin{align*}
\left(\begin{array}{cc} a & b\\ -b & a\end{array}\right)\in\GL_2(A).
\end{align*}
Diagonalising, $\SSS_{\CC}$ is isomorphic to $\GG^2_m$, with complex conjugation on $\SSS(\CC)$ corresponding to $(z_1,z_2)\mapsto(\overline{z_2},\overline{z_1})$. Therefore, the elements of $\SSS(\RR)$ map to the elements $(z,\overline{z})$, stable under conjugation. More generally, the characters of $\SSS_{\CC}$ are the homomorphisms 
\begin{align*}
(z_1,z_2)\mapsto z^p_1z^q_2,
\end{align*} 
for any $(p,q)\in\ZZ\times\ZZ$, with complex conjugation acting as $(p,q)\mapsto (q,p)$.

Consequently, to give a representation of $\SSS$ on a real vector space $V$ is the same as to give a $\ZZ\times\ZZ$-grading of $V(\CC)$ such that $\overline{V^{p,q}}=V^{q,p}$ for all $p$ and $q$, which is precisely the definition of a Hodge structure on $V$. We thus define morphisms, tensor products and duals of Hodge structures as morphisms, tensor products and duals of representations of $\SSS$. We normalise the relation so that $(z_1,z_2)$ acts on $V^{p,q}$ as $z^{-p}_1z^{-q}_2$. A complex structure on a real vector space $V$ is then precisely a Hodge structure $\SSS\rightarrow\GL(V)$ coming from a homomorphism $\CC\rightarrow\End(V)$. 

For $n\in\ZZ$ and $R=\ZZ$, $\QQ$ or $\RR$, we let $R(n)$ be the ($R$-)Hodge structure $V=R$, where $\SSS$ acts on $V(\RR)=\RR$ by the character $(z\overline{z})^n$ and, hence, 
\begin{align*}
V(\CC)=V_{-n}(\CC).
\end{align*}
This is referred to as a {\it Tate twist}. For an ($R$-)Hodge structure $V$ of weight $n$, a {\it Hodge tensor} is a multilinear form $t:V^r\rightarrow R$ such that the map
\begin{align*}
V\otimes V\otimes\cdots\otimes V\rightarrow R(-nr/2)
\end{align*}
is a morphism of Hodge structures. 

If we denote by $C:=h(i)$ the {\it Weil operator}, then a {\it polarisation} on $V$ is a Hodge tensor
\begin{align*}
\psi:V\times V\rightarrow R
\end{align*}
such that
\begin{align*}
\psi_C:V(\RR)\times V(\RR)\rightarrow\RR:(x,y)\mapsto\psi(x,Cy) 
\end{align*}
is symmetric and positive definite. A {\it polarisation} on an ($R$-)Hodge structure $V=\oplus_n V_n$ is a system $(\psi_n)_n$ of polarisations on the $V_n$.

\section{Abelian varieties}

Consider an {\it Abelian variety} $A$ over $\CC$ of dimension $g$. Then $A$ is isomorphic to a {\it complex torus} $\CC^g/\Lambda$, where $\Lambda$ is the $\ZZ$-module generated by an $\RR$-basis for $\CC^g$. The isomorphism $\Lambda\otimes\RR\cong\CC^g$ defines a complex structure on $\Lambda\otimes\RR$ and there exists an alternating form 
\begin{align*}
\psi:\Lambda\times\Lambda\rightarrow\ZZ
\end{align*} 
such that $\psi_{\RR}(x,Cy)$ is symmetric and positive definite and 
\begin{align*}
\psi_{\RR}(Cx,Cy)=\psi_{\RR}(x,y),
\end{align*} 
for all $x,y\in\Lambda\otimes\RR$. In other words, $\Lambda\cong H_1(A,\ZZ)$ is a $\ZZ$-Hodge structure of weight $-1$ equipped with a polarisation. In fact, by \cite{M}, Theorem 6.8, the functor $A\mapsto H_1(A,\ZZ)$ is an equivalence from the category of Abelian varieties over $\CC$ to the category of polarised $\ZZ$-Hodge structures of type $(-1,0),(0,-1)$. Therefore, the answer to the question of the previous section is that {\it one can study the problem of parameterising Abelian varieties in terms of Hodge structures}. 

Consider the case of Abelian varieties of dimension one, otherwise known as {\it elliptic curves}. An elliptic curve over $\CC$ is the quotient of $\CC$ by a free $\ZZ$-module $\Lambda$ of rank $2$. Two elliptic curves $\CC/\Lambda$ and $\CC/\Lambda'$ are isomorphic if and only if $\Lambda'=\alpha\Lambda$ for some $\alpha\in\CC^{\times}$. We summarise the perspective explained in \cite{H}:

Often, when considering elliptic curves, we fix $\CC$ and vary $\Lambda$. Instead, however, we may fix $\Lambda:=\ZZ^2$ and vary the complex structure on $\ZZ^2\otimes\RR=\RR^2$ i.e. we vary the morphism
\begin{align*}
h:\CC^{\times}\rightarrow\GL_2(\RR)
\end{align*}
extending to a homomorphism $\CC\rightarrow{\rm M}_2(\RR)$ of $\RR$-algebras. Given such a morphism, we obtain an isomorphism of complex vector spaces $i_h:\RR^2\rightarrow\CC$ defined by 
\begin{align*}
i^{-1}_h(z)=h(z)\cdot i^{-1}_h(1):=h(z)\cdot e_0,
\end{align*}
where we choose $e_0=(1,0)\in\RR^2$. The quotient $\CC/i_h(\ZZ^2)$ is an elliptic curve.

Therefore, let
\begin{align*}
h_0:\CC^{\times}\rightarrow\GL_2(\RR):a+ib\mapsto\left(\begin{array}{cc} a & b\\ -b & a\end{array}\right)
\end{align*}
and let $h:=\gamma h_0\gamma^{-1}$, where 
\begin{align*}
\gamma=\left(\begin{array}{cc} x & y\\ w & z\end{array}\right)\in\GL_2(\RR)^+.
\end{align*}
Note that, for any such $h$, the standard symplectic form given by
\begin{align*}
(u,v)\mapsto u^t\left(\begin{array}{cc} 0 & -1\\ 1 & 0\end{array}\right)v
\end{align*}
is a polaristion for the corresponding $\ZZ$-Hodge structure.

For $h_0(z)$, the $z$-eigenspace in $\RR^2\otimes\CC$ is the complex subspace generated by $(-i,1)$. The $\overline{z}$-eigenspace is its complex conjugate, generated by $(i,1)$. Therefore, for $h(z)$, the $z$-eigenspace is generated by
\begin{align*}
\left(\begin{array}{cc} x & y\\ w & z\end{array}\right)\left(\begin{array}{c}-i\\ 1\end{array}\right)=\left(\begin{array}{c} -xi+y\\ -wi+z\end{array}\right)
\end{align*}
or, equivalently, $(\overline{\tau}_h,1)$, where $\tau_h:=xi+y/wi+z$, and the $\overline{z}$-eigenspace is generated by $(\tau_h,1)$. Note that this latter subspace is precisely the middle term in the filtration associated to the $\ZZ$-Hodge structure given by $h$.

Now, $i_h$ extends $\CC$-linearly to a map
\begin{align*}
i_{h,\CC}:\RR^2\otimes\CC=\CC\cdot\left(\begin{array}{c} \overline{\tau}_h \\ 1\end{array}\right)\oplus\CC\cdot\left(\begin{array}{c} \tau_h\\ 1\end{array}\right)\rightarrow\CC
\end{align*}
and, since it commutes with the action of $\CC$ on both sides, we deduce that $i_{h,\CC}$ is the quotient of $\RR^2\otimes\CC$ by the $\overline{z}$-eigenspace. Therefore, since $i_h(e_0)=1$ and $i_h((0,1))=i_h(-\tau_h e_0+(\tau_h,1))=-\tau_h$, 
\begin{align*}
i_h(\ZZ^2)=\ZZ\oplus\ZZ\tau_h.
\end{align*}

We conclude that $\CC/i_h(\ZZ^2)$ varies over all isomorphism classes of elliptic curves as $h$ varies over the $\GL_2(\RR)^+$-conjugacy class of $h_0$. The map $h\mapsto\tau_h$ is a $\GL_2(\RR)^+$-equivariant bijection between this conjugacy class and $\HH$.

For Abelian varieties of dimension $g$, the situation is similar. We replace $\ZZ^2$ by $\ZZ^{2g}$ and fix the standard symplectic form given by
\begin{align*}
-J:=\left(\begin{array}{cc} 0 & -{\rm id}\\ {\rm id} & 0\end{array}\right).
\end{align*}
We let
\begin{align*}
h_0:\CC^{\times}\rightarrow\GL_{2g}(\RR):a+bi\mapsto a+bJ,
\end{align*}
which factors through the group 
\begin{align*}
\GSp_{2g}(\RR)=\{g\in\GL_{2g}(\RR):g^tJg=\nu(g) J\},
\end{align*}
where $\nu:\GSp_{2g}\rightarrow\GG_m$ is a homomorphism of linear algebraic groups.
The $\GSp_{2g}(\RR)^+$-conjugacy class of $h_0$ corresponds to the set of $\ZZ$-Hodge structures on $\ZZ^{2g}$ having type $(-1,0),(0,-1)$ for which $J$ induces a polarisation. Using the description of the Hodge filtration, as in the case of elliptic curves, one can identify this set in a $\GSp_{2g}(\RR)^+$-equivariant manner with a Hermitian symmetric domain 
\begin{align*}
\HH_g:=\{Z=X+iY\in M_{g\times g}(\CC):Z=Z^t, Y>0\}
\end{align*} 
called the {\it Siegel upper half-space of genus $g$}. 

\section{The Siegel upper half-space}

Let us return then to our account of Hodge structures. Having fixed a $g\in\NN$, we denote the Hodge structure corresponding to a point $\tau\in\HH_g$ by $V_{\tau}$ and we denote the corresponding Hodge filtration by $F_{\tau}$. For any given $(p,q)\in\ZZ\times\ZZ$, the dimension $d(p,q)$ of $V^{p,q}_{\tau}$ is constant as $\tau$ varies over $\HH_g$ and we have a continuous map
\begin{align*}
\tau\mapsto [V^{p,q}_{\tau}]:\HH_g\rightarrow {\rm G}_{d(p,q)}(V(\CC)),
\end{align*}
from $\HH_g$ to the complex, projective variety of $d(p,q)$-dimensional subspaces of $V(\CC)$. 

The subspace dimensions of $F_{\tau}$ are then also constant as $\tau$ varies over $\HH_g$ and, if we denote by ${\rm F}_d(V(\CC))$ the complex, projective variety parameterising such filtrations of $V(\CC)$, then the map
\begin{align*}
f:\tau\mapsto [F_{\tau}]:\HH_g\rightarrow {\rm F}_d(V(\CC))
\end{align*}
is holomorphic. In light of these properties, we refer to the set of Hodge structures corresponding to the points of $\HH_g$ as a {\it holomorphic family of Hodge structures}. 

Finally, the differential of $f$ at $\tau$ is a $\CC$-linear map
\begin{align*}
df_{\tau}:T_{\tau}\HH_g\rightarrow T_{[F_{\tau}]}{\rm F}_d(V(\CC))
\end{align*}
from the tangent plane of $\HH_g$ at $\tau$ to the tangent plane of ${\rm F}_d(V(\CC))$ at $[F_{\tau}]$.
By \cite{M}, (17), $T_{[F_{\tau}]}{\rm F}_d(V(\CC))$ is a subset of
\begin{align*}
\bigoplus_p\Hom(F^p_{\tau},V(\CC)/F^p_{\tau})
\end{align*}
but, in this case, the image of $df_{\tau}$ is actually contained in the space
\begin{align*}
\bigoplus_p\Hom(F^p_{\tau},F^{p-1}_{\tau}/F^p_{\tau})
\end{align*} 
and we say that this holomorphic family of Hodge structures is a {\it variation of Hodge structures}.

\section{Families of Hodge structures}

The above situation can be abstracted as follows: let $V$ be a finite dimensional $\RR$-vector space and let $T$ be a finite set of tensors on $V$, including a nondegenerate bilinear form $t_0$. Fix an $n\in\NN$ and let 
\begin{align*}
d:\ZZ\times\ZZ\rightarrow\NN
\end{align*} 
be a symmetric function such that $d(p,q)=0$ for almost all $(p,q)$, including every $(p,q)$ such that $p+q\neq n$. 

Consider the set $S(d,T)$ of Hodge structures on $V$ such that, for all $(p,q)\in\ZZ\times\ZZ$,
\begin{align*}
\dim V^{p,q}=d(p,q),
\end{align*} 
every $t\in T$ is a Hodge tensor and $t_0$ is a polarisation. This is naturally a subspace of 
\begin{align*}
\prod_{(p,q):d(p,q)\neq 0}{\rm G}_{d(p,q)}(V(\CC)).
\end{align*} 
Therefore, $S(d,T)$ can be given the subspace topology and, by \cite{M}, Theorem 2.14, (assuming it is non-empty) any connected component has a unique complex structure such that the corresponding set of Hodge structures constitute a holomorphic family. Furthermore, if such a family is actually a variation of Hodge structures, then the corresponding connected component $S^+$ has the structure of a Hermitian symmetric domain. In fact, every Hermitian symmetric domain is of the form $S^+$ for a suitable $V$, $T$ and $d$.

\section{The algebraic group}

Recall the topological space $S(d,T)$ from the previous section and let $S^+$ be a connected component. Fix a point $h_0\in S^+$ and let $G$ be the smallest algebraic subgroup of $\GL(V)$ such that 
\begin{align*}
h:\SSS\rightarrow\GL(V)
\end{align*} 
factors through $G$ for every $h\in S^+$ i.e. the intersection of all subgroups having this property. As in the proof of \cite{M}, Theorem 2.14 (a), for any $g\in G(\RR)^+$, $gh_0g^{-1}\in S^+$ and, in fact, the map
\begin{align*}
g\mapsto gh_0g^{-1}:G(\RR)^+\rightarrow S^+
\end{align*}
is surjective. In other words, $S^+$ is the $G(\RR)^+$-conjugacy class of $h_0$.

\section{Shimura data}

Motivated by our example of Abelian varieties, we want to consider $\ZZ$-(or $\QQ$)-Hodge structures. This will be achieved by choosing an algebraic group $G$ defined over $\QQ$ and embedding this into $\GL(V)$ for some $\QQ$-vector space $V$. The $\ZZ$-structure will come from the choice of a lattice in $V$. 

\begin{defini}
A Shimura datum is a pair $(G,X)$, where $G$ is a reductive group over $\QQ$ and $X$ is a $G(\RR)$-conjugacy class of morphisms $h:\SSS\rightarrow G_{\RR}$ such that, for one (or, equivalently, all) $h\in X$,
\begin{itemize}
\item Only the characters $z\mapsto 1$, $z\mapsto z/\bar{z}$ and $z\mapsto \bar{z}/z$ occur in the representation of $\SSS$ on the Lie algebra of $G^{\ad}_{\CC}$.
\item Conjugation by $h(i)$ is a Cartan involution of $G^{\ad}$.
\item For every simple factor $H$ of $G^{\ad}$, the map $\SSS\rightarrow H_{\RR}$ is not trivial.
\end{itemize} 
\end{defini}
By a {\it reductive} algebraic group we refer to a connected, linear algebraic group with trivial unipotent radical. The {\it unipotent radical} of a linear algebraic group is the unipotent part of its radical, where its {\it radical} is the neutral component of its maximal normal, solvable subgroup. The semisimple groups are those linear algebraic groups with trivial radical. In particular, they are reductive.

Now let $(G,X)$ be a Shimura datum. By the first of the axioms above, $\GG_{m}(\RR)=\RR^{\times}$, which is naturally a subgroup of $\SSS(\RR)=\CC^{\times}$, acts trivially on $\mathfrak{g}_{\CC}$. As the action of $G$ on $\mathfrak{g}$ factors through $G^{\ad}$ and the action of $G^{\ad}$ is faithful, the image of $\RR^{\times}$ in $G(\RR)$ must belong to the centre. In particular, the restriction of any $h\in X$ to $\GG_m$ is independent of $h$ and we refer to its reciprocal $w$ as the {\it weight homomorphism} since, for any representation $\rho:G_{\RR}\rightarrow\GL(V)$, $\rho\circ w$ defines the weight decomposition of the Hodge structure given by $\rho\circ h$ on $V$. 

Now let $\rho:G_{\RR}\rightarrow\GL(V)$ be a faithful representation. By \cite{M}, Proposition 5.9, $X$ has a unique structure of a complex manifold such that the family of Hodge structures induced on $V$ by $\rho\circ h$ as $h$ varies over $X$ is holomorphic. In fact, the first axiom implies that it is a variation of Hodge structures. Therefore, from our earlier discussion of families of Hodge structures, $X$ is a finite disjoint union of Hermitian symmetric domains.

Alternatively, consider a connected component $X^+$ of $X$. By \cite{M}, Proposition 5.7 (a), we may consider $X^+$ as a $G^{\ad}(\RR)^+$-conjugacy class of morphisms $\SSS\rightarrow G^{\ad}_{\RR}$. Let $h\in X^+$ and decompose $G^{\ad}_{\RR}$ into a product of simple factors $H_i$ so that $h=(h_i)_i$, where $h_i$ is the projection of $h$ to $H_i$. By \cite{M}, Lemma 4.7, if $H_i(\RR)$ is compact then $h_i$ is trivial. Otherwise, given the conditions satisfied by $h$, there exists a Hermitian symmetric domain $D_i$ such that $H_i(\RR)^+$ coincides with $\Hol(D_i)^+$ and $D_i$ is in natural one-to-one correspondence with the $H_i(\RR)^+$-conjugacy class $X^+_i$ of $h_i$. Therefore, the product $D$ of the $D_i$ is a Hermitian symmetric domain on which $G^{\ad}(\RR)^+$ acts via a surjective homomorphism $G^{\ad}(\RR)^+\rightarrow\Hol(D)^+$ with compact kernel and there is a natural identification of $D$ with $X^+=\prod_i X^+_i$.
 
\begin{defini}
A morphism of Shimura data 
\begin{align*}
(G_1,X_1)\rightarrow (G_2,X_2)
\end{align*} 
is a morphism $\phi:G_1\rightarrow G_2$ such that, for every $h\in X_1$, $\phi\circ h\in X_2$. If $\phi$ is a closed immersion, we refer to $(G_1,X_1)$ as a Shimura subdatum.
\end{defini}

\begin{defini}
Let $(G,X)$ be a Shimura datum. Let $X^{\ad}$ be the $G^{\ad}(\RR)$-conjugacy class of morphisms $\SSS\rightarrow G^{\ad}_{\RR}$ containing the image of $X$. Then $(G^{\ad},X^{\ad})$ is a Shimura datum called the adjoint Shimura datum and 
\begin{align*}
(G,X)\rightarrow(G^{\ad},X^{\ad})
\end{align*} 
is a morphism of Shimura data. 
\end{defini}

\section{Congruence subgroups}
Let $G$ be a reductive subgroup of $\GL_n$ defined over $\QQ$. We denote by $G(\ZZ)$ the group $G(\QQ)\cap\GL_n(\ZZ)$. Recall the following definition, independent of the embedding of $G$ in $\GL_n$:

\begin{defini}
A subgroup $\Gamma$ of $G(\QQ)$ is arithmetic if $\Gamma\cap G(\ZZ)$ has finite index in $\Gamma$ and $G(\ZZ)$ i.e. if $\Gamma$ and $G(\ZZ)$ are commensurable.
\end{defini}

Now suppose that $(G,X)$ is a Shimura datum. We would like to consider the corresponding Hodge structures up to isomorphism and this is the role of the group $\Gamma$. We may also wish to distinguish additional structure to that already encoded in the group $G$. The most obvious such structure is distinguished by the following class of arithmetic subgroups:

\begin{defini}
The principal congruence subgroup of level $N$ is defined as the group
\begin{align*}
\Gamma(N):=\{g\in G(\ZZ): g\equiv{\rm id}\text{ mod }N\},
\end{align*}
where the congruence relation is entry-wise. 
\end{defini}

In the case of Abelian varieties, where $G=\GSp_{2g}$ and we consider the $\ZZ$-Hodge structure on $\Lambda=H_1(A,\ZZ)$, the group $\Gamma(N)$ also distinguishes between different bases for the {\it $N$-torsion subgroup} $\frac{1}{N}\Lambda/\Lambda$, rather than simply the isomorphism class of $\Lambda$ along with its polarisation.

Of course, the definition of the principal congruence subgroup depends on the embedding of $G$ in $\GL_n$. Therefore, we define a {\it congruence subgroup} of $G(\QQ)$ to be a subgroup containing some $\Gamma(N)$ as a subgroup of finite index. This notion does not depend on the embedding.

\section{Adeles}

The ring of {\it finite (rational) ad\`eles} $\AAA_f$ comprises the elements 
\begin{align*}
\alpha=(\alpha_p)\in\prod_p\QQ_p
\end{align*} 
such that, for almost all primes $p$, $\alpha_p\in\ZZ_p$. It is endowed with the topology for which a basis of open sets are those of the form $\prod_pU_p$, where $U_p$ is open in $\QQ_p$, and $U_p=\ZZ_p$ for almost all $p$. Similarly, for an algebraic group $G$, defined over $\QQ$, one can choose an embedding into $\GL_n$ and define $G(\AAA_f)$ as those elements 
\begin{align*}
g=(g_p)_p\in\prod_p G(\QQ_p)
\end{align*} 
such that $g_p\in\GL_n(\ZZ_p)$ for almost all $p$. However, this definition of $G(\AAA_f)$ is independent of the embedding into $\GL_n$ and so is the basis of open sets, defined analogously to the above. 

By \cite{M}, Proposition 4.1, for any compact open subgroup $K$ of $G(\AAA_f)$, $K\cap G(\QQ)$ is a congruence subgroup $\Gamma$ of $G(\QQ)$ and every congruence subgroup arises this way. Loosely speaking, considering the congruence relation defining $\Gamma$ prime-by-prime gives rise to $K$ and vice-versa.

Later, we will also need the more general definition of $\AAA_{E,f}$, the {\it finite ad\`eles over a number field} $E$, which we define as $\AAA_f\otimes E$ or, equivalently, as the ring of elements 
\begin{align*}
\alpha=(\alpha_{\upsilon})\in\prod_{\upsilon}E_{\upsilon},
\end{align*} 
over all finite places $\upsilon$ of $E$ such that, for almost all $\upsilon$, $\alpha_{\upsilon}\in\OOO_{E_\upsilon}$. The {\it ad\`ele ring} $\AAA_E$ arises when we include factors for the infinite places of $E$. Therefore, any $\alpha\in\AAA_E$ can be written as a pair $(\alpha_{\infty},\alpha_f)$, where $\alpha_f\in\AAA_{E,f}$.

\section{Neatness}

Let $G$ be an algebraic subgroup of $\GL_n$ defined over $\QQ$. The following definition is independent of the embedding into $\GL_n$: 

\begin{defini}
An element $g\in G(\QQ)$ is neat if the subgroup of $\overline{\QQ}^{\times}$ generated by its eigenvalues is torsion free. 
\end{defini}

One says that a congruence subgroup $\Gamma$ is {\it neat} if all of its elements are neat. There is also a notion of neatness for compact open subgroups of $G(\AAA_f)$, for which we refer the reader to \cite{KY}, 4.1.4. In particular, if $K$ is neat then so is the congruence subgroup $G(\QQ)\cap gKg^{-1}$, for any $g\in G(\AAA_f)$. Every compact open subgroup $K$ of $G(\AAA_f)$ contains a neat compact open subgroup $K'$ with finite index.

\section{Shimura varieties}
Finally, we give the definition of a Shimura variety:
\begin{defini}
Let $(G,X)$ be a Shimura datum and let $K$ be a compact open subgroup of $G(\AAA_f)$. The Shimura variety attached to $(G,X)$ and $K$ is the double coset space
\begin{align*}
\Sh_K(G,X)(\CC):=G(\QQ)\backslash X\times (G(\AAA_f)/K).
\end{align*}
\end{defini}

This definition invariably seems abstruse at first. However, it is a simple calculation to see that
\begin{align*}
\Sh_K(G,X)(\CC)=\coprod_{g\in\cal{C}}\Gamma'_g\backslash X,
\end{align*}
where $\cal{C}$ is a set of representatives for the double coset space $G(\QQ)\backslash G(\AAA_f)/K$ and $\Gamma'_g:=G(\QQ)\cap gKg^{-1}$ is a congruence subgroup. Note that, by \cite{PR}, Theorem 5.1, $\cal{C}$ is a finite set. However, since we are interested in connected components, choose a connected component $X^+$ of $X$ and denote by $G(\QQ)_+$ its stabiliser in $G(\QQ)$. Then
\begin{align*}
\Sh_K(G,X)(\CC)=\coprod_{g\in\cal{C}_+}\Gamma_g\backslash X^+,
\end{align*}
where $\cal{C}_+$ is a set of representatives for the double coset space $G(\QQ)_+\backslash G(\AAA_f)/K$ and $\Gamma_g:=G(\QQ)_+\cap gKg^{-1}$. By \cite{M}, Lemma 5.12, $\cal{C}_+$ is also a finite set.

\section{Complex structure}

Any arithmetic subgroup $\Gamma$ of $G(\QQ)$ acts on $X$ through $G^{\ad}(\QQ)$ and, by \cite{M}, Proposition 3.2, its image is also arithmetic. For any arithmetic subgroup $\Gamma$ of $G(\QQ)$, the intersection $\Gamma\cap G(\QQ)_+$ acts on $X^+$. We say that its image under the map $G^{\ad}(\RR)^+\rightarrow\Hol(X^+)^+$ is an {\it arithmetic subgroup} of $\Hol(X^+)^+$. 

If $\Gamma$ is neat then the image of $\Gamma\cap G(\QQ)_+$ in $\Hol(X^+)^+$ is neat and, in particular, torsion free. By \cite{M}, Proposition 3.1, such an arithmetic subgroup of $\Hol(X^+)^+$ acts freely on $X^+$ and the corresponding quotient has a unique complex structure such that the quotient map is a local isomorphism. In general then, $\Gamma\backslash X^+$ has the structure of a (possibly singular) complex analytic variety.

\section{Algebraic structure}

The fundamental result of Baily and Borel \cite{BB} states that the quotient of $X^+$ by any torsion free, arithmetic subgroup of $\Hol(X^+)^+$ has a canonical realisation as a complex, quasi-projective, algebraic variety. In particular, if $K$ is neat, $\Sh_K(G,X)(\CC)$ is the analytification of a quasi-projective variety $\Sh_K(G,X)_{\CC}$. 

A further theorem of Borel \cite{B} states that, for any smooth, quasi-projective variety $V$ over $\CC$, any holomorphic map from $V(\CC)$ to $\Sh_K(G,X)(\CC)$ is regular. For example, given any inclusion $K_1\subset K_2$ of neat compact open subgroups of $G(\AAA_f)$, we have a natrual morphism of algebraic varieties
\begin{align*}
\Sh_{K_1}(G,X)_{\CC}\rightarrow\Sh_{K_2}(G,X)_{\CC}.
\end{align*}
Therefore, varying $K$, we get an inverse system of algebraic varieties 
\begin{align*}
(\Sh_K(G,X)_{\CC})_K
\end{align*} 
and we write the scheme-theoretic limit of this system as $\Sh(G,X)_{\CC}$. On the system there is a natural action of $G(\AAA_f)$ given by
\begin{align*}
\cdot g:\Sh_K(G,X)(\CC)\rightarrow\Sh_{g^{-1}Kg}(G,X)(\CC):[x,a]_K\mapsto[x,ag]_{g^{-1}Kg},
\end{align*}
where we use $[\cdot,\cdot]_K$ to denote a double coset belonging to $\Sh_K(G,X)(\CC)$. By the theorem of Borel, this action is regular on components. Therefore, for any given $g\in G(\AAA_f)$, we obtain an algebraic correspondence
\begin{align*}
\Sh_K(G,X)_{\CC}\leftarrow\Sh_{K\cap gKg^{-1}}(G,X)_{\CC}\xrightarrow{\cdot g}\Sh_{g^{-1}Kg\cap K}(G,X)_{\CC}\rightarrow\Sh_{K}(G,X)_{\CC},
\end{align*}
where the outer maps are the natural projections. We refer to this correspondence as a {\it Hecke correspondence}.

Finally, if we have a morphism
\begin{align*}
f:(G_1,X_1)\rightarrow (G_2,X_2)
\end{align*} 
of Shimura data and two compact open subgroups $K_1\subset G_1(\AAA_f)$ and $K_2\subset G_2(\AAA_f)$ such that $f(K_1)\subset K_2$, then we obtain a morphism
\begin{align*}
\Sh_{K_1}(G_1,X_1)(\CC)\rightarrow\Sh_{K_2}(G_2,X_2)(\CC),
\end{align*}
which, again by the theorem of Borel, is a regular map
\begin{align*}
\Sh_{K_1}(G_1,X_1)_{\CC}\rightarrow\Sh_{K_2}(G_2,X_2)_{\CC}.
\end{align*}
We refer to the images of such maps as {\it Shimura subvarieties}. We also have an induced morphism
\begin{align*}
\Sh(G_1,X_1)_{\CC}\rightarrow\Sh(G_2,X_2)_{\CC}
\end{align*}
of the limits, by which we mean an inverse system of regular maps, compatible with the actions of $G_1(\AAA_f)$ and $G_2(\AAA_f)$.

\section{Special subvarieties}

{\it Special subvarieties} constitute the smallest class of irreducible algebraic subvarieties containing the connected components of Shimura subvarieties and closed under taking irreducible components of images under Hecke correspondences. The precise definition is the following:

\begin{defini}
Let $\Sh_K(G,X)_{\CC}$ be a Shimura variety. A closed irreducible subvariety $Z$ is called special if there exists a morphism of Shimura data 
\begin{align*}
(G',X')\rightarrow(G,X)
\end{align*} 
and $g\in G(\AAA_f)$ such that $Z$ is an irreducible component of the image of
\begin{align*}
\Sh(G',X')_{\CC}\rightarrow\Sh(G,X)_{\CC}\xrightarrow{\cdot g}\Sh(G,X)_{\CC}\rightarrow\Sh_K(G,X)_{\CC}.
\end{align*}

\end{defini}

The situation is analogous to the case of Abelian varieties, where the special subvarieties are the Abelian subvarieties and their translates under torsion points.

By definition, if we let $K'\subset G(\AAA_f)$ be a compact open subgroup contained in $K$ and consider the natural morphism of Shimura varieties
\begin{align*}
\pi:\Sh_{K'}(G,X)_{\CC}\rightarrow\Sh_{K}(G,X)_{\CC},
\end{align*}
\begin{itemize}
\item if $Z$ is a special subvariety of $\Sh_{K'}(G,X)_{\CC}$, then $\pi(Z)$ is a special subvariety of $\Sh_K(G,X)_{\CC}$.
\item if $Z$ is a special subvariety of $\Sh_{K}(G,X)_{\CC}$, then any irreducible component of $\pi^{-1}Z$ is a special subvariety of $\Sh_{K'}(G,X)_{\CC}$.
\end{itemize} 

\section{Special points}

The natural definition of a {\it special point} is then the following:

\begin{defini}
A special point in $\Sh_K(G,X)_{\CC}$ is a special subvariety of dimension zero.
\end{defini}
However, we can characterise special points in a more concrete manner: consider a special point $[h,g]_K\in\Sh_K(G,X)(\CC)$. Let $M:=\MT(h)$ be the {\it Mumford-Tate group} of $h$ i.e. the smallest algebraic subgroup $H$ of $G$ (defined over $\QQ$) such that $h:\SSS\rightarrow G_{\RR}$ factors through $H_{\RR}$ and let $X_M$ denote the orbit $M(\RR)\cdot h$ inside $X$. Then $(M,X_M)$ is a Shimura subdatum of $(G,X)$ and, if we let $X^+_M$ be the connected component $M(\RR)^+\cdot h$ of $X_M$, then the image of $X^+_M\times\{g\}$ in $\Sh_K(G,X)(\CC)$ defines the smallest special subvariety containing $[h,g]_K$. Therefore, $X_M$ must be zero dimensional and so $M$ must be commutative. It is a general fact that any subgroup of $G$ defined over $\QQ$ and containing $h(\SSS)$ is reductive. Therefore, $M$ is a torus. 

On the other hand if $T$ is a torus in $G$ and $h\in X$ factors through $T_{\RR}$ then $[h,g]_K\in\Sh_K(G,X)(\CC)$ is clearly a special point for any $g\in G(\AAA_f)$. Therefore, we may define a special point as any point $[h,g]_K\in\Sh_K(G,X)(\CC)$ such that $\MT(h)$ is a torus. Of course, the choice of $h$ is only well-defined up to conjugation by an element of $G(\QQ)$, but this doesn't affect the property of $\MT(h)$ being a torus.

\section{Canonical model}

It is possible to define a model for $\Sh_K(G,X)_{\CC}$ that is canonical in a sense one can make precise. As we have seen, $\Sh_K(G,X)(\CC)$ is often a moduli space for Abelian varieties and the main theorem of complex multiplication gives us a description of how Galois groups act on sets of CM-Abelian varieties. Therefore, we would like the Galois action on $\Sh_K(G,X)(\CC)$ to agree with this description, whenever it applies. In order to achieve this, the canonical model satisfies a generalised version of this description given in terms of Deligne's group-theoretic $(G,X)$ language. We provide a very brief summary of the theory explained more thoroughly in \cite{M}, \S12, \S13 and \S14. 

Recall that a {\it model} over a number field $E$ for a complex algebraic variety $V$ is a variety $V_0$ defined over $E$ with an isomorphism $\phi:V_{0,\CC}\rightarrow V$, though we will follow convention and omit any mention of this isomorphism. First we define the field of definition $E:=E(G,X)$ of the canonical model. It is referred to as the {\it reflex field} and, as we will see, it does not depend on $K$. This independence is one reason for having several connected components in the definition of a Shimura variety. 

For a subfield $k$ of $\CC$, we write $\mathcal{C}(k)$ for the set of $G(k)$-conjugacy classes of cocharacters of $G_k$ defined over $k$ i.e.
\begin{align*}
\mathcal{C}(k)=G(k)\backslash\Hom(\GG_{m,k},G_k).
\end{align*}
Any homomorphism $k\rightarrow k'$ induces a map $\mathcal{C}(k)\rightarrow\mathcal{C}(k')$, so $\Aut(k'/k)$ acts on $\mathcal{C}(k')$. 

For $h\in X$, we obtain a cocharacter
\begin{align*}
\mu_h:\GG_{m,\CC}\xrightarrow{z\mapsto(z,1)}\GG^2_{m,\CC}\cong\SSS_{\CC}\xrightarrow{h_{\CC}}G_{\CC}
\end{align*}
of $G_{\CC}$ and so the $G(\RR)$-conjugacy class $X$ of $h$ maps to an element $c(X)\in\mathcal{C}(\CC)$. The reflex field $E$ is then the fixed field of the stabiliser of $c(X)$ in $\Aut(\CC)$. By what follows, we will see that $E$ is a number field.

Suppose that 
\begin{align*}
[h,g]_K\in\Sh_K(G,X)(\CC) 
\end{align*} 
is a special point i.e. $M:=\MT(h)$ is a torus. Therefore, since all cocharacters of $M$ are defined over $\overline{\QQ}$ and $\mu_h$ factors through $M_{\CC}$, $\mu_h$ is defined over a finite extension $E_h$ of $\QQ$. Note that $E_h$ does not depend on the choice of $h$. By \cite{M}, Remark 12.3 (b), $E$ is contained in $E_h$.

For any $t\in M(E_h)$, the element 
\begin{align*}
\prod_{\sigma:E_h\rightarrow\overline{\QQ}}\sigma(t)
\end{align*}
is stable under $\Gal(\overline{\QQ}/\QQ)$ and so belongs to $M(\QQ)$. The so-called {\it reciprocity morphism} is defined by
\begin{align*}
r_h:\AAA^{\times}_{E_h,f}\rightarrow M(\AAA_f):a\mapsto\prod_{\sigma:E_h\rightarrow\overline{\QQ}}\sigma(\mu_h(a)).
\end{align*}
Finally, recall the (surjective) {\it Artin map}
\begin{align*}
{\rm Art}_{E_h}:\AAA^{\times}_{E_h}\rightarrow\Gal(E_h^{\ab}/E_h)
\end{align*}
from class field theory and let ${\rm Art}^{-1}_{E_h}$ denote its reciprocal.

\begin{defini}
We say that a model of $\Sh_K(G,X)_{\CC}$ over $E$ is canonical if every special point $[h,g]_K$ in $\Sh_K(G,X)(\CC)$ has coordinates in $E_h^{\ab}$ and
\begin{align*}
\sigma[h,g]_K=[h,r_h(s_f)a]_K,
\end{align*} 
for any $\sigma\in\Gal(E_h^{\ab}/E_h)$ and $s=(s_{\infty},s_f)\in\AAA^{\times}_{E_h}$ such that ${\rm Art}^{-1}_{E_h}(s)=\sigma$.
\end{defini}
By \cite{M}, Theorem 13.7, if a canonical model exists, it is unique up to unique isomorphism. The difficult theorem is that canonical models actually exist. For a discussion, see \cite{M}, \S14. 

A {\it model} of $\Sh(G,X)_{\CC}$ over $E$ is an inverse system of varieties over $E$, endowed with a right action of $G(\AAA_f)$, which over $\CC$ is isomorphic to $\Sh(G,X)_{\CC}$ with its $G(\AAA_f)$ action. Such a system is {\it canonical} if each component is canonical in the above previous sense. 

By \cite{M}, Theorem 13.7 (b), if for all compact open subgroups $K$ of $G(\AAA_f)$ $\Sh_K(G,X)_{\CC}$ has a canonical model, then so does $\Sh(G,X)_{\CC}$ and it is unique up to unique isomorphism. In particular, by \cite{M}, Theorem 13.6, the action of $G(\AAA_f)$ is defined over $E$. By \cite{M}, Remark 13.8, if $(G',X')\rightarrow (G,X)$ is a morphism of Shimura data and $\Sh(G',X')_{\CC}$ and $\Sh(G,X)_{\CC}$ have canonical models, then the induced morphism
\begin{align*}
\Sh(G',X')_{\CC}\rightarrow\Sh(G,X)_{\CC}
\end{align*}
is defined over $E(G',X')\cdot E(G,X)$.

\section{The Andr\'e-Oort conjecture}

The Andr\'e-Oort conjecture is the following statement regarding the geometry of Shimura varieties:

\begin{conj}\label{AO}
Let $(G,X)$ be a Shimura datum, $K$ a compact open subgroup of $G(\AAA_f)$ and $\Sigma$ a set of special points in $\Sh_K(G,X)(\CC)$. Then every irreducible component of the Zariski closure of $\cup_{s\in\Sigma}s$ in $\Sh_K(G,X)_{\CC}$ is a special subvariety.
\end{conj}

In the remainder of this article, we are going to apply the Pila-Zannier strategy to the Andr\'e-Oort conjecture. The Andr\'e-Oort conjecture is analogous to the Manin-Mumford conjecture (first proved by Raynaud \cite{R}), asserting that the irreducible components of the Zariski closure of a set of torsion points in an Abelian variety are the translates of Abelian subvarieties by torsion points. The task at hand is essentially to combine a number of different ingredients. We follow the outline given by Ullmo in \cite{U}, \S5 for the case of $\mathcal{A}^r_6$.

\section{Reductions}

Let $Y$ denote an irreducible component of the Zariski closure of $\cup_{s\in\Sigma}s$ in $\Sh_K(G,X)_{\CC}$. Let $[h,g]_K\in Y$ denote a point such that $M:=\MT(h)$ is maximal among such groups. Note that the maximality is independent of the choice of $h$. We say that such a point is {\it Hodge generic} in $Y$. 

Let $X_M:=M(\RR)\cdot h$. Then, by \cite{EY}, Proposition 2.1, $Y$ is contained in the image of the morphisms
\begin{align*}
\Sh_{K_M}(M,X_M)_{\CC}\rightarrow\Sh_{gKg^{-1}}(G,X)_{\CC}\xrightarrow{\cdot g}\Sh_K(G,X)_{\CC},
\end{align*}
where $K_M:=M(\AAA_f)\cap gKg^{-1}$. Denote by $f$ their composition and let $Y_M$ be an irreducible component of $f^{-1}Y$. Then $Y$ is a special subvariety of $\Sh_K(G,X)_{\CC}$ if and only if $Y_M$ is a special subvariety of $\Sh_{K_M}(M,X_M)_{\CC}$. Furthermore, $Y_M$ is Hodge generic in $\Sh_{K_M}(M,X_M)_{\CC}$. Therefore, we may assume that $Y$ is Hodge generic in $\Sh_K(G,X)_{\CC}$.

Let $(G^{\ad},X^{\ad})$ be the adjoint Shimura datum associated to $(G,X)$ and let $K^{\ad}$ be a compact open subgroup of $G^{\ad}(\AAA_f)$ containing the image of $K$. Then $Y$ is a special subvariety of $\Sh_K(G,X)_{\CC}$ if and only if its image $Y^{\ad}$ in $\Sh_{K^{\ad}}(G^{\ad},X^{\ad})_{\CC}$ is a special subvariety. Furthermore, if $Y$ is Hodge generic in $\Sh_K(G,X)_{\CC}$, then $Y^{\ad}$ is Hodge generic in $\Sh_{K^{\ad}}(G^{\ad},X^{\ad})_{\CC}$. Therefore, we may assume that $G$ is semisimple of adjoint type.

Recall that the irreducible components of the image of a special subvariety under a Hecke correspondence are again special subvarieties. Therefore, if we fix a connected component $X^+$ of $X$, we may assume that $Y$ is contained in the image $S:=\Gamma\backslash X^+$ of $X^+\times\{1\}$ in $\Sh_K(G,X)(\CC)$, where $\Gamma:=G(\QQ)_+\cap K$. We denote a point in $S$ as $[h]$ for some $h\in X^+$. 

\section{Galois orbits}

The first ingredient is a lower bound for the size of the Galois orbit of a special point. By the definition of special subvarieties, the choice of $K$ is irrelevant in the Andr\'e-Oort conjecture. Thus, we may assume that $K$ is neat and a product of compact open subgroups $K_p$ in $G(\QQ_p)$. 

Now let $[h]\in S$ be a special point. Recall that $M:=\MT(h)$ is a torus and let $L$ denote its {\it splitting field}, by which we mean the smallest field over which $M$ becomes isomorphic to a product of the multiplicative group. Note that this is a finite, Galois extension of $\QQ$ containing $E_h$ and is independent of the choice of $h$. 

Let $K_M$ denote the compact open subgroup $M(\AAA_f)\cap K$ of $M(\AAA_f)$, which is equal to the product of the $K_{M,p}:=M(\QQ_p)\cap K_p$. Let $K^m_M$ be the maximal compact open subgroup of $M(\AAA_f)$, which is unique since $M$ is a torus and equal to the product of the maximal compact open subgroups $K^m_{M,p}$ of $M(\QQ_p)$. Note that $K_{M,p}=K^m_{M,p}$ for almost all primes $p$. The following conjecture is a natural generalisation of \cite{EMO}, Problem 14, posed by Edixhoven for $\mathcal{A}_g$:
\begin{conj}\label{gal}
There exist positive constants $c_1$, $B_1$ and $\mu_1$ such that, for any special point $[h]\in S$, 
\begin{align*}
|\Gal(\overline{\QQ}/L)\cdot[h]|>c_1B_1^{i(M)}[K^m_M:K_M]D^{\mu_1}_L,
\end{align*}
where $i(M)$ is the number of places such that $K_{M,p}\neq K^m_{M,p}$ and $D_L$ is the absolute value of the discriminant of $L$.
\end{conj} 
Note that, although the groups $K^m_H$ and $K_M$ depend on the choice of $h$, they are well-defined up to conjugation by an element of $\Gamma$ and, hence, the index $[K^m_M:K_M]$ is well-defined. By \cite{UY2}, Th\'eor\`eme 6.1, this bound is known to hold under the generalised Riemann hypothesis for CM fields and, by \cite{T}, Theorem 1.1, it holds unconditionally in the case of $\mathcal{A}_g$, for $g$ at most $6$.

\section{Realisations}

We refer to a point $h\in X^+$ as a {\it pre-special point} if $[h]\in S$ is a special point. The second ingredient in the Pila-Zannier strategy is an upper bound for the height of a pre-special point in a fundamental domain $\mathcal{F}$ of $X^+$ with respect to $\Gamma$. As opposed to the case of an Abelian variety, this is a non-trivial issue. 

For a sensible notion of height, we must first choose a {\it realisation} $\mathcal{X}$ of $X^+$. By this we mean an analytic subset of a complex, quasi-projective variety $\mathcal{\wt{X}}$, with a transitive holomorphic action of $G(\RR)^+$ on $\mathcal{X}$ such that, for any $x_0\in\mathcal{X}$, the orbit map \begin{align*}
G(\RR)^+\rightarrow\mathcal{X}:g\mapsto g\cdot x_0
\end{align*} 
is semi-algebraic and identifies $\mathcal{X}$ with $G(\RR)^+/K_{\infty}$, where $K_{\infty}$ is a maximal compact subgroup of $G(\RR)^+$ (recall that $G$ is semisimple and adjoint). A morphism of realisations is then a $G(\RR)^+$-equivariant  biholomorphism. By \cite{U}, Lemme 2.1, any realisation has a canonical semi-algebraic structure and any morphism of realisations is semi-algebraic. Therefore, $X^+$ has a canonical semi-algebraic structure.

A subset $Z\subset\mathcal{X}$ is called an {\it irreducible algebraic subvariety} of $\mathcal{X}$ if $Z$ is an irreducible component of the analytic set $\mathcal{X}\cap\wt{Z}$, where $\wt{Z}$ is an algebraic subset of $\mathcal{\wt{X}}$. By \cite{U}, Lemme 2.1, $\mathcal{X}\cap\wt{Z}$ has finitely many analytic components and they are semi-algebraic. 
Also note that, by \cite{KUY}, Corollary B.1, this notion is independent of our choice of $\mathcal{X}$. In particular, we have a well defined notion of an {\it irreducible algebraic subvariety} of $X^+$.

\section{Heights}

For the remainder of this article, we will fix as our realisation the so-called {\it Borel embedding} of $X^+$ into its {\it compact dual} $X^{\vee}$. We refer to \cite{UY3}, 3.3 for the following definitions:

As before, for a point $h\in X^+$, let
\begin{align*}
\mu_h:\GG_{m,\CC}\xrightarrow{z\mapsto(z,1)}\GG^2_{m,\CC}\cong\SSS_{\CC}\xrightarrow{h_{\CC}} G_{\CC}
\end{align*}
be the corresponding cocharacter and let $M_X$ be the $G(\CC)$-conjugacy class of $\mu_h$. Let $V$ be a faithful representation of $G$ on a finite dimensional $\QQ$-vector space so that, for each point $h\in X^+$, we obtain a Hodge structure $V_h$ and a Hodge filtration
\begin{align*}
F_h:=\{\cdots\supset F^p_h\supset F^{p+1}_h\supset\cdots\},\ F^p_h:=\oplus_{r\geq p}V^{r,s}_h.
\end{align*}

Fix a point $h_0\in X^+$ and let $P$ be the {\it parabolic subgroup} of $G(\CC)$ stabilising $F_{h_0}$. We define $X^{\vee}$ to be the complex, projective variety $G(\CC)/P$, which is naturally a subvariety of the flag variety $\Theta_{\CC}:=\GL(V_{\CC})/Q$, where $Q$ is the parabolic subgroup of $\GL(V_{\CC})$ stabilising $F_{h_0}$. Therefore, we have a surjective map from $M_X$ to $X^{\vee}$ sending $\mu_h$ to $F_h$.

The Borel embedding $X\hookrightarrow X^{\vee}$ is the map $h\mapsto F_h$. It is injective since, by \cite{M}, \S2, (18), the Hodge filtration determines the Hodge decomposition. In other words, the maximal compact subgroup $K_{\infty}$ of $G(\RR)^+$ constituting the stabiliser of $h_0$ is equal to $G(\RR)^+\cap P$.

However, $\Theta_{\CC}$ has a natural model $\Theta$ over $\QQ$ such that, for any extension $L$ of $\QQ$, a point of $\Theta(L)$ corresponds to a filtration defined over $L$. By definition, $X^{\vee}$ is defined over the reflex field $E:=E(G,X)$ and, by the proof of \cite{UY3}, Proposition 3.7, a special point $h\in X^+$ is defined over the splitting field of a maximal torus $T$ of $\GL(V)$ such that $T_{\CC}$ contains the Mumford-Tate group of $h$.  

Therefore, since a pre-special point $h\in X^+$ has algebraic coordinates, we are allowed to talk about its {\it (multiplicative) height} $H(h)$, as defined in \cite{BG}, Definition 1.5.4. The following is a natural generalisation of \cite{PT2}, Theorem 3.1, due to Tsimerman:
\begin{conj}\label{height}
There exist positive constants $c_2$, $B_2$, $\mu_2$ and $\mu_3$ such that, for any pre-special point $h\in\mathcal{F}$,
\begin{align*}
H(h)<c_2B_2^{i(M)}[K^m_M:K_M]^{\mu_2}D^{\mu_3}_L.
\end{align*} 
\end{conj}

Finally, let $h\in X^+$ be a pre-special point and let $L$ be the splitting field of a maximal torus $T$ of $\GL(V)$ such that $T_{\CC}$ contains the Mumford-Tate group of $h$. The dimension $d$ of $T$ is at most the dimension of $V$ and the Galois action on the character group of $T$ is given by a homomorphism
\begin{align*}
\Gal(L/\QQ)\hookrightarrow\GL_d(\ZZ).
\end{align*}
Since, by a classical result of Minkowski, the number of isomorphism classes of finite groups contained in $\GL_d(\ZZ)$ is finite, the degree of $L$ is bounded by a positive constant depending only on $G$.

\section{Definability}

In order to apply the Pila-Wilkie counting theorem, one requires the following theorem: 
\begin{teo}\label{def}
The restriction $\pi_{|\mathcal{F}}$ of the uniformisation map 
\begin{align*}
\pi:X^+\rightarrow S
\end{align*} 
is definable in $\RR_{\rm an,exp}$.
\end{teo}
This theorem is discussed in several articles. It was first proved for restricted theta functions by Peterzil and Starchenko \cite{PS}. In particular, this addressed the case of $\mathcal{A}_g$. It is known for general Shimura varieties due to the work of Klingler, Ullmo and Yafaev \cite{KUY}. 

\section{Ax-Lindemann-Weierstrass}

The final ingredient is the hyperbolic Ax-Lindemann-Weierstrass conjecture. In order to state the conjecture, we require the notion of a {\it weakly special subvariety}:

\begin{defini}
A variety $V$ in $S$ is weakly special if the (analytic) connected components of $\pi^{-1}V$ are algebraic in $X^+$.
\end{defini}

This definition is actually the characterisation \cite{UY3}, Theorem 1.2 of the original definition \cite{UY3}, Definition 2.1. However, given some familiarity with Shimura varieties, the proof is fairly straightforward and this characterisation is precisely what we need. The term {\it weakly special} is motivated by the fact that all special subvarieties are weakly special whereas, as explained in \cite{Mo}, weakly special subvarieties are special subvarieties if and only if they contain a special point. 

\begin{teo}\label{AL}
Let $Z$ be an algebraic subvariety of $S$. Maximal, irreducible, algebraic subvarieties of $\pi^{-1}Z$ are precisely the irreducible components of the preimages of maximal, weakly special subvarieties contained in $Z$.
\end{teo}

Again, this problem and its history are discussed at length in several other articles. The theorem above is due to Klingler, Ullmo and Yafaev \cite{KUY}. It was first proven for compact Shimura varieties by Ullmo and Yafaev \cite{UY} and for $\mathcal{A}_g$ by Pila and Tsimerman \cite{PT}.

\section{Pila-Wilkie}

Let $A\subset\RR^m$ be a definable set in an o-minimal structure and let $A^{\rm alg}$ be the union of all connected, positive dimensional, semi-algebraic subsets contained in $A$. Recall the Pila-Wilkie counting theorem, first proved for rational points in \cite{PW} and later for algebraic points in \cite{P2}: 

\begin{teo}\label{PilaWilkie}
For every $\epsilon>0$ and $k\in\NN$, there exists a positive constant $c$, depending only on $A$, $k$ and $\epsilon$, such that, for any real number $T\geq 1$, the number of points lying on $A\setminus A^{\rm alg}$, whose coordinates in $\RR^m$ are algebraic of degree at most $k$ and of multiplicative height at most $T$, is at most $cT^{\epsilon}$.
\end{teo} 
In this article, the o-minimal structure will be $\RR_{\rm an,exp}$ and {\it definable} will always mean definable in $\RR_{\rm an,exp}$. 

\section{Final reduction}

The final reduction is the following result due to Ullmo, appearing as Theorem 4.1 in \cite{U}:

\begin{teo}
Let $Z$ be a Hodge generic subvariety of $\Sh_K(G,X)_{\CC}$, strictly contained in $S$. Suppose that, if $S$ is a product $S_1\times S_2$ of connected components of Shimura varieties, then $Z$ is not of the form $S_1\times Z'$, for a subvariety $Z'$ of $S_2$. Then the union of all positive-dimensional, weakly special subvarieties of $\Sh_K(G,X)_{\CC}$ contained in $Z$ is not Zariski dense in $Z$. 
\end{teo}

We apply the theorem to $Y$ noting that the assumption in the theorem is no loss of generality: if necessary, we simply replace $S$ by $S_2$ and $Y$ by $Y'$. Thus, we may assume that the union of all positive-dimensional special subvarieties of $\Sh_K(G,X)_{\CC}$ contained in $Y$ is not Zariski dense in $Y$. 

Therefore, if we are able to show that all but a finite number of special points in $Y$ lie on a positive-dimensional special subvariety of $\Sh_K(G,X)_{\CC}$ contained in $Y$, then the theorem implies that $Y=S$.

\section{The Pila-Zannier strategy}

By Theorem \ref{def}, $\pi_{|\mathcal{F}}$ is definable and so
\begin{align*}
\wt{Y}:=\pi^{-1}Y\cap\mathcal{F}
\end{align*}
is a definable set. By assumption, $Y$ contains a dense set of special points and so is defined over a finite extension $F$ of $E$.

Consider a pre-special point $h\in\wt{Y}$ and let $L$ denote the splitting field of $M:=\MT(h)$. The Galois orbit $\Gal(\overline{\QQ}/LF)\cdot[h]$ is contained in $Y$ and, if Conjecture \ref{gal} holds, then
\begin{align*}
|\Gal(\overline{\QQ}/LF)\cdot [h]|>c'_1B_1^{i(M)}[K^m_M:K_M]D^{\mu_1}_L,
\end{align*}
where $c'_1:=c_1/[F:E]$. On the other hand, by \cite{M}, Example 12.4 (a),  $\Gal(\overline{\QQ}/LF)\cdot[h]$ is contained in the image of the morphism
\begin{align*}
\Sh_{K_M}(M,h)(\CC)\rightarrow\Sh_K(G,X)(\CC),
\end{align*}
induced by the inclusion of Shimura data. Therefore, let 
\begin{align*}
[h,m]_K\in\Sh_K(G,X)(\CC)
\end{align*} 
denote an element of $\Gal(\overline{\QQ}/LF)\cdot[h]$, where $m\in M(\AAA_f)$ is given by the explicit description of the Galois action. Since $[h,m]_K\in S$, $m$ is equal to $qk$, for some $q\in G(\QQ)$ and $k\in K$. Denote by $h'$ the point of $\wt{Y}$ such that $[h']=[h,m]_K$. Then, up to conjugation by an element of $\Gamma$, 
\begin{align*}
M':=\MT(q^{-1}\cdot h)=q^{-1}Mq
\end{align*}
is equal to $\MT(h')$ and
\begin{align*}
K^m_{M'}/K_{M'}=q^{-1}K^m_Mq/q^{-1}M(\AAA_f)q\cap K.
\end{align*} 
Conjugation by $q$ yields a bijection between this quotient and
\begin{align*}
K^m_M/M(\AAA_f)\cap qKq^{-1},
\end{align*}
which has cardinality $[K^m_M:K_M]$ since $q=mk^{-1}$.

Consequently, if Conjecture \ref{height} holds, then
\begin{align*}
H(h')<c_2B_2^{i(M)}[K^m_M:K_M]^{\mu_2}D^{\mu_3}_L.
\end{align*} 
Therefore, since all pre-special points in $X^+$ have algebraic co-ordinates of bounded degree, Theorem \ref{PilaWilkie} implies that, for any $\epsilon>0$, there exists a constant $c$, depending only on $\wt{Y}$ and $\epsilon$, such that there are at most 
\begin{align*}
c(c_2B_2^{i(M)}[K^m_M:K_M]^{\mu_2}D^{\mu_3}_L)^{\epsilon}
\end{align*} 
pre-special points on $\wt{Y}\setminus\wt{Y}^{\rm alg}$ belonging to $\Gal(\overline{\QQ}/LF)\cdot[h]$. 

Therefore, we may choose $\epsilon$ sufficiently small such that, if either $[K^m_M:K_M]$ or $D_L$ is large enough, then there exists a point in $\Gal(\overline{\QQ}/LF)\cdot[h]$ such that the corresponding point $h'\in\wt{Y}$ belongs to a positive dimensional, semi-algebraic set contained in $\wt{Y}$. Therefore, by \cite{KUY}, Lemma B.2, $h'$ belongs to an irreducible algebraic subvariety of $X^+$ contained in $\wt{Y}$ and so, by Theorem \ref{AL} (the hyperbolic Ax-Lindemann-Weierstrass theorem), there exists a weakly special subvariety $V$ contained in $Y$ such that $[h']\in V$. Therefore, $V$ is a special subvariety of positive dimension and $[h]$ belongs to a special subvariety contained in $Y$.

Therefore, on Y, in the complement of all positive dimensional, special subvarieties contained in $Y$, the quantities $[K^m_M:K_M]$ and $D_L$ corresponding to special points are bounded. By \cite{UY4}, Proposition 3.21, the set of tori equal to the Mumford-Tate group of a pre-special point such that $[K^m_M:K_M]$ and $D_L$ are bounded lie in only finitely many $\Gamma$-conjugacy classes. In particular, such pre-special points lie above only finitely points in $S$.

\end{document}